\input amstex
\documentstyle{amsppt}
\topmatter
\magnification=\magstep1
\pagewidth{5.2 in}
\pageheight{6.7 in}
\abovedisplayskip=10pt
\belowdisplayskip=10pt
\NoBlackBoxes
\title
 $q$-Euler and Genocchi numbers
\endtitle
\author  Taekyun Kim \endauthor
\affil\rm{{Institute of Science Education,}\\
{ Kongju National University, Kongju 314-701, S. Korea}\\
{e-mail: tkim$\@$kongju.ac.kr (or tkim64$\@$hanmail.net)}}\\
\endaffil

\abstract{Carlitz has introduced an interesting $q$-analogue of
Frobenius-Euler numbers in [4]. He has indicated a corresponding
Stadudt-Clausen theorem and also some interesting congruence
properties of the $q$-Euler numbers. In this paper we give another
construction of $q$-Euler numbers, which are different than his
$q$-Euler numbers. By using our $q$-Euler numbers, we define the
$q$-analogue of Genocchi numbers and investigate the relations
between $q$-Euler numbers and $q$-analogs of Genocchi numbers. }
\endabstract
\thanks 2000 Mathematics Subject Classification  11S80, 11B68 \endthanks
\thanks Key words and phrases:  Sums of powers, Bernoulli number, Bernoulli polynomials \endthanks
\rightheadtext{ Taekyun Kim    } \leftheadtext{ $q$-Euler and
Genocchi numbers }
\endtopmatter

\document

\head 1. Introduction \endhead Throughout this paper, we consider
a complex number $q\in\Bbb C$ with $|q|<1 $ as an indeterminate.
The $q$-analogue of $n$ is defined by $[n]_q=\frac{1-q^n}{1-q}.$
The ordinary Euler numbers are defined by the generating function
as follows:
$$F(t)=\frac{2}{e^t+1}=e^{Et}=\sum_{n=0}^{\infty}E_n\frac{t^n}{n!},
\text{ $|t|<\pi$,}\tag1$$ where we use the technique method
notation by replacing $E^m$ by $E_m$ ($m\geq0$), symbolically,
cf.[2, 6].

From Eq.(1), we can derive the Genocchi numbers as follows:
$$G(t)=\frac{2t}{e^t +1}=\sum_{n=0}^{\infty} G_n \frac{t^n}{n!},
\text{ $|t|<\pi$ .}\tag2$$ It satisfies $G_1=1,$
$G_3=G_5=G_7=\cdots=0,$ and even coefficients are given
$G_m=2(1-2^{2m})B_{2m}=2mE_{2m-1},$ where $B_m$ are the $m$-th
ordinary Bernoulli numbers, cf.[6]. It follows from (2) and
Stadudt-Clasusen  theorem that Genocchi numbers are integers. For
$x\in \Bbb R$(=the field of real numbers) the Euler polynomials
are defined by
$$F(x,t)=F(t)e^{xt}=\frac{2}{e^t+1}e^{xt}=\sum_{n=0}^{\infty}E_{n}(x)\frac{t^n}{n!},
\text{  ($|t|<\pi$)}. \tag3$$ From (3), we can also derive the
definition of Genocchi polynomials as follows:
$$\frac{2t}{e^t+1}e^{xt}=\sum_{n=0}^{\infty}G_n(x)\frac{t^n}{n!},
\text{ ($|t|<\pi$). }\tag4$$ The following formulae ((5)-(6)) are
well known in [6].
$$E_m(x)=\sum_{k=0}^m\binom mk \frac{G_{k+1}}{k+1}x^{m-k}. \tag 5$$
For $n, m\geq 1 ,$ and $n$ odd, we have
$$(n^m-n)G_m=\sum_{k=1}^{m-1}\binom mk n^k G_k Z_{m-k}(n-1)
,\tag6$$ where $Z_m(n)=1^m-2^m+3^m-\cdots+(-1)^{n+1}n^m .$ In this
paper we give the $q$-analogs of the above Eq.(5) and Eq.(6). The
purpose of this paper is to give  another construction of
$q$-Euler numbers, which are different than a $q$-Eulerian numbers
of Carlitz.  From the definition of our $q$-Euler numbers, we
derive the $q$-analogs of Genocchi numbers and investigate the
properties of $q$-Genocchi numbers which are related to $q$-Euler
numbers.

\head 2. $q$-Euler Numbers and polynomials
\endhead

Let $q$ be a complex number with $q<1$. In [3, 4] Carlitz
constructed $q$-analogue of Eulerian numbers. We now consider
another construction of a $q$-Eulerian numbers, which are
different than his $q$-Eulerian numbers. First we consider the
following generating functions:
$$F_q(t)=[2]_qe^{\frac{t}{1-q}}\sum_{j=0}^{\infty}\frac{(-1)^j}{1+q^{j+1}}\left(\frac{1}{1-q}\right)^j\frac{t^j}{j!}
=e^{E_qt}=\sum_{n=0}^{\infty}E_{n,q}\frac{t^n}{n!}, \tag7$$ and
$$F_{q}(x,t)=[2]_qe^{\frac{t}{1-q}}\sum_{j=0}^{\infty}\frac{(-1)^jq^{jx}}{1+q^{j+1}}\left(\frac{1}{1-q}\right)^j\frac{t^j}{j!}
=e^{E_q(x)t}=\sum_{n=0}^{\infty}E_{n,q}(x)\frac{t^n}{n!}, \tag8$$
where we use the technique method notation by replacing $E_q^n$ by
$E_{n,q} ,$ symbolically. Thus we have
$$\aligned
&E_{n,q}=\frac{[2]_q}{(1-q)^n}\sum_{l=0}^n\binom
nl\frac{(-1)^l}{1+q^{l+1}}, \\
&E_{n,q}(x)=\frac{[2]_q}{(1-q)^n}\sum_{l=0}^n\binom
nl\frac{(-1)^l}{1+q^{l+1}}q^{lx},\endaligned\tag 8-1$$ where
$\binom nl$ is binomial coefficient.

 By (8-1), we easily see that $\lim_{q\rightarrow 1}E_{n,q}=E_n$
and $\lim_{q\rightarrow 1}E_{n,q}(x)=E_n(x) .$ From Eq.(8), we can
derive the below Eq.(9):
$$F_q(x,t)=[2]_q\sum_{n=0}^{\infty}(-1)^nq^ne^{[n+x]_qt}=\sum_{n=0}^{\infty}E_{n,q}(x)\frac{t^n}{n!}.\tag9$$
By (9), we easily see that
$$
E_{n,q}(x)=\frac{[2]_q}{[2]_{q^m}}[m]_q^n\sum_{a=0}^{m-1}(-1)^aq^aE_{n,q^m}(\frac{a+x}{m})=\sum_{k=0}^n\binom
nk [x]_q^{n-k}q^{kx}E_{k,q} , \text{  for $m$ odd}.\tag10$$ This
is equivalent to
$$[2]_{q^m}E_{n,q}(xm)=[2]_q[m]_q^n\sum_{a=0}^{m-1}(-1)^aq^aE_{n,q^m}(\frac{a}{m}+x),
\text{ for $m$ odd}. \tag11$$ If we put $x=0$ in Eq.(11), then we
have
$$[m]_{-q}E_{n,q}-[m]_q^n\frac{[m(n+1)]_{-q}}{[n+1]_{-q}}E_{n,q^m}
=\sum_{l=0}^{n-1}\binom nl
[m]_q^lE_{l,q^m}\sum_{a=1}^{m-1}(-1)^aq^{a(l+1)}[a]_q^{n-l},
\tag12$$ where $[m]_{-q}=\frac{1+q^m}{1+q}$ for $m$ odd.

Define the operation $*$ on $f_n(q)$ as follows:
$$(1-[m]_q^n)*f_n(q)=[m]_{-q}f_n(q)-[m]_q^n\frac{[m(n+1)]_{-q}}{[n+1]_{-q}}f_n(q^m).
\tag13$$ By (12) and (13), we obtain the following:
 \proclaim{ Proposition 1} For $m,n\in \Bbb N$ and $m$ odd, we have
 $$(1-[m]_q^n)*E_{n,q}=\sum_{l=0}^{n-1}\binom nl [m]_q^l
 E_{l,q^m}\sum_{a=1}^{m-1}(-1)^aq^{a(l+1)}[a]_q^{n-l}. $$
\endproclaim
For any positive integer $n$, it is easy to see that
$$-[2]_q\sum_{l=0}^{\infty}(-1)^{l+n}q^{l+n}e^{[l+n]_qt}+[2]_q\sum_{l=0}^{\infty}(-1)^lq^le^{[l]_qt}
=[2]_q\sum_{l=0}^{n-1}(-1)^lq^le^{[l]_qt}.\tag14$$ From (9) and
(14) we can derive the below:
$$\sum_{l=0}^{n-1}(-1)^lq^l[l]_q^m=\frac{1}{[2]_q}\left((-1)^{n+1}q^nE_{m,q}(n)-E_{m,q}\right).$$
Therefore we obtain the following: \proclaim{Proposition 2} For
$n, m \in\Bbb N$, we have
$$\sum_{l=0}^{n-1}(-1)^lq^l[l]_q^m=\frac{1}{[2]_q}\left((-1)^{n+1}q^nE_{m,q}(n)-E_{m,q}\right).$$
\endproclaim
In the recent many authors have studies the sums of powers of
consecutive integers, cf.[1, 5, 7, 10, 11]. The above Proposition
2 is the another $q$-analogue of the sums of powers of consecutive
integers. The Genocchi numbers $G_n$ are defined by the generating
function:
$$G(t)=\frac{2t}{e^t+1}=e^{Gt}=\sum_{n=0}^{\infty}G_n\frac{t^n}{n!},
\text{ ($|t|<\pi$) ,} $$ where we use the technique method
notation by replacing $G^m$ by $G_m$ ($m\geq 0$), symbolically. It
satisfies $G_1=1 ,$ $G_3=G_5=G_7=\cdots=0$ and even coefficients
are given $G_m=2(1-2^{2m})B_{2m}=2mE_{2m-1},$ cf.[6]. We now
derive the $q$-extension of the above Genocchi numbers from the
definition of our $q$-Euler numbers.

\head 3. $q$-Genocchi Numbers and polynomials
\endhead
By the meaning of (1) and (2), let us define the $q$-extension of
Genocchi numbers as follows:
 $$G_q(t)=[2]_qt\sum_{n=0}^{\infty}(-1)^nq^ne^{[n]_qt}=\sum_{n=0}^{\infty}G_{n,q}\frac{t^n}{n!},
 \text{ ($|t|<\pi$) }.$$
Note that $\lim_{q\rightarrow
1}G_q(t)=\frac{2t}{e^t+1}=\sum_{n=0}^{\infty}G_n\frac{t^n}{n!}.$
Hence, $\lim_{q\rightarrow 1}G_{n,q}=G_n. $ In [8], the
$q$-Bernoulli numbers are defined by
$$-t\sum_{n=0}^{\infty}q^n
e^{[n]_qt}=\sum_{n=0}^{\infty}B_{n,q}\frac{t^n}{n!}, \text{ (
$|t|<2\pi$) .} \tag 15$$ It was known that $\lim_{q\rightarrow
1}B_{n,q}=B_n ,$ cf.[8, 9]. By (15), we easily see that
$$-[2]_qt\sum_{n=0}^{\infty}q^ne^{[n]_qt}+2[2]_qt\sum_{n=0}^{\infty}q^{2n}e^{[2n]_qt}
=[2]_qt\sum_{n=0}^{\infty}(-1)^nq^ne^{[n]_qt}.\tag16$$ From (15)
and (16), we can derive the below Eq.(17):
$$G_{n,q}=[2]_qB_{n,q}-2[2]_q^nB_{n,q^2}. \tag17$$
Let us consider the $q$-analogue of Genocchi polynomials as
follows:
$$G_q(x,t)=[2]_qt\sum_{n=0}^{\infty}(-1)^nq^{n+x}e^{[n+x]_qt}=\sum_{n=0}^{\infty}G_{n,q}(x)
\frac{t^n}{n!}. \tag18$$ By (18), we easily see that
$$G_q(x,t)=[2]_qq^{xt}e^{\frac{t}{1-q}}\sum_{l=0}^{\infty}\frac{(-1)^l}{1+q^{l+1}}q^{lx}\left(\frac{1}{1-q}
\right)^l\frac{t^l}{l!}.\tag19$$ Thus, we have
$$G_{n,q}(x)=n\left(\frac{1}{1-q}\right)^{n-1}\sum_{l=0}^{n-1}\binom{n-1}{l}\frac{(-1)^l}{1+q^{l+1}}q^{(l+1)x}.
$$ From (8) and (19), we can derive the below equality:
$$\aligned
F_q(x,t)&=[2]_q\sum_{n=0}^{\infty}(-1)^nq^ne^{[n+x]_qt}
=\frac{e^{[x]_qt}}{q^xt}[2]_qq^x
t\sum_{n=0}^{\infty}(-1)^nq^ne^{q^x[n]_qt}\\
&=e^{[x]_qt}\sum_{n=0}^{\infty}q^{nx}\frac{G_{n+1,q}}{n+1}\frac{t^n}{n!}
=\sum_{n=0}^{\infty}\left(\sum_{k=0}^n\binom nk
[x]_q^{n-k}q^{nx}\frac{G_{n+1,q}}{n+1}\right)\frac{t^n}{n!}.
\endaligned\tag20$$
By (20), we easily see that
$$E_{n,q}(x)=\sum_{k=0}^n\binom nk
[x]_q^{n-k}q^{nx}\frac{G_{n+1,q}}{n+1}. \tag21$$ Remark. The
Eq.(21) is the $q$-analogue of Eq.(5).

Therefore we obtain the following theorem:
 \proclaim{Theorem 3}
For any positive integer $n$, we have
$$\aligned
&(a)  \text{ }
G_{n,q}(x)=n\left(\frac{1}{1-q}\right)^{n-1}\sum_{l=0}^{n-1}\binom{n-1}l
\frac{(-1)^l}{1+q^{l+1}}q^{(l+1)x},\\
& (b)  \text{ } E_{n,q}(x)=\sum_{k=0}^n\binom nk
[x]_q^{n-k}q^{nx}\frac{G_{n+1,q}}{n+1}, \\
& (c) \text{ } G_{n,q}=[2]_qB_{n,q}-2[2]_q^nB_{n,q^2},
\endaligned$$
where $B_{n,q}$ are the $q$-Bernoulli numbers which are defined in
[8].\endproclaim By (18), we easily see that
$$\sum_{n=0}^{\infty}G_{n,q}(x)\frac{t^n}{n!}=\sum_{n=0}^{\infty}
\left(\frac{[2]_q}{[2]_{q^m}}[m]_q^{n-1}\sum_{a=0}^{m-1}(-1)^aq^{a+x}G_{n,q^m}(\frac{x+a}{m})\right)\frac{t^n}{n!},
\text{ for $m$ odd }.\tag22$$ Thus we obtain the following:
\proclaim{Theorem 4} Let $m\in\Bbb N$ and $m$ odd. Then we see
that
$$G_{n,q}(x)=\frac{[2]_q}{[2]_{q^m}}[m]_q^{n-1}\sum_{a=0}^{m-1}(-1)^aq^{a+x}G_{n,q^m}(\frac{x+a}{m})
=\sum_{k=0}^{\infty}\binom nk q^{kx}G_{k,q}[x]_q^{n-k}. $$ This is
equivalent to
$$G_{n,q}(mx)=\frac{[2]_q}{[2]_{q^m}}[m]_q^{n-1}\sum_{a=0}^{m-1}(-1)^aq^{a+mx}
G_{n,q^m}(x+\frac{a}{m}). \tag23$$
\endproclaim
If we take $x=0$ in Eq.(23), then we easily see that
$$[2]_{q^m}[m]_qG_{n,q}-[2]_q[m]_q^nG_{n,q^m}\frac{[2]_{q^{m(n+1)}}}{[2]_{q^{n+1}}}
=[2]_q\sum_{k=0}^{n-1} \binom nk
[m]_q^kG_{k,q^m}\sum_{a=0}^{m-1}(-1)^aq^{a(k+1)}[a]_q^{n-k}.\tag24$$
From the definition of the operation $*$ in the previous section,
we note that
$$\left([m]_q-[m]_q^n\right)*f_n(q)=[2]_{q^m}[m]_qf_n(q)-[2]_q[m]_q^n\frac{[2]_{q^{m(n+1)}}}
{[2]_{q^{n+1}}}f_n(q^m).\tag25$$ By (24) and (25), we easily see
that
$$\left([m]_q-[m]_q^n\right)*G_{n,q}=[2]_q\sum_{k=0}^{n-1}\binom
nk [m]_q^kG_{n,q^m}\sum_{a=0}^{m-1}(-1)^aq^{a(k+1)}[a]_q^{n-k}. $$

\enddocument